\pdfoutput=1
\RequirePackage{ifpdf}
\ifpdf 
\documentclass[pdftex]{sigma}
\else
\documentclass{sigma}
\fi

\numberwithin{equation}{section}
\graphicspath{{./Images2/}}

\newtheorem{Theorem}{Theorem}[section]
\newtheorem{Corollary}[Theorem]{Corollary}
\newtheorem{Lemma}[Theorem]{Lemma}
 { \theoremstyle{definition}
\newtheorem{Definition}[Theorem]{Definition}
\newtheorem{Example}[Theorem]{Example}
\newtheorem{Remark}[Theorem]{Remark} }

\begin{document}
\allowdisplaybreaks

\newcommand{\arXivNumber}{1810.05777}

\renewcommand{\PaperNumber}{119}

\FirstPageHeading

\ShortArticleName{Counting Collisions in an $N$-Billiard System Using Angles Between Collision Subspaces}

\ArticleName{Counting Collisions in an $\boldsymbol{N}$-Billiard System\\ Using Angles Between Collision Subspaces}

\Author{Sean GASIOREK}

\AuthorNameForHeading{S.~Gasiorek}

\Address{School of Mathematics and Statistics, Carslaw F07, University of Sydney, NSW 2011, Australia}
\Email{\href{mailto:sean.gasiorek@sydney.edu.au}{sean.gasiorek@sydney.edu.au}}
\URLaddress{\url{https://www.seangasiorek.com}}

\ArticleDates{Received July 24, 2020, in final form November 10, 2020; Published online November 21, 2020}

\Abstract{The principal angles between binary collision subspaces in an $N$-billiard system in $d$-dimensional Euclidean space are computed. These angles are computed for equal masses and arbitrary masses. We then provide a bound on the number of collisions in the planar 3-billiard system problem. Comparison of this result with known billiard collision bounds in lower dimensions is discussed. }


\Classification{37J37; 70F99; 55R80; 70F16}

\section{Introduction}	

Mathematical billiards is a well-studied example of a dynamical system, and its geometric and dynamic properties are of interest to mathematicians and physicists alike \cite{KT, Tab}. Boltzmann first described the kinetic theory of a low-density gas as a system of interactions of small particles
(e.g., atoms, molecules), which can simply be modelled as a billiard system. In particular, Boltzmann studied the many-particle problem by considering subsets of particles of size 2,~3,~\dots whose motion was not affected by particles outside of these subsets. Boltzmann specifically studied the case when these subsets each included a single binary collision. The natural next step is the subcollection of three particles or spheres. This is addressed in \cite{MC1,MC2}, and it is shown that the maximum number of collisions amongst three hard spheres in $\mathbb{R}^{d}$ is four. The maximum number of collisions is unknown in the case of four or more spheres in any dimension higher than 1.

Finding estimates for the maximum number of collisions of a billiard system is highly dependent upon properties of the given system: the quantity of billiard balls, the underlying space, their respective masses and radii all affect the total number of collisions. For example, the configuration space of two equal-radii billiard balls on the same side of a fixed wall in 1-dimensional space is isomorphic to the motion of a single point-billiard inside a wedge of angle measure~$\alpha < \pi$~\cite{Gal2003,Tab}. In such a setting, the maximum number of collisions of the single billiard in~the~wedge is $\lceil \pi/\alpha \rceil$. Higher-dimensional analogues for estimates on the number of collisions in a~polyhedral angle similarly depend upon geometric properties of the bounding hyperplanes, see, e.g., \cite{Sev,Sin}. Further, alternate formulations of billiards (e.g.,~\cite{AT2018,BMT2020}) may produce different bounds on the number of collisions in a multi-billiard system.

In \cite{BFK1,BFK2} a uniform bound for the number of collisions in semi-dispersing billiards in terms of the minimum and maximum masses and radii of a collection of $N$ billiard balls is computed. We provide an alternate approach to bounding the number of collisions in the $N$-billiard system introduced in \cite{FKM}, focusing on the planar $N=3$ case.

This paper is constructed as follows. Section \ref{S2} outlines the construction of $N$-body billiards, an alternate construction of the configuration space, and defines the geometric properties of angles linear subspaces of a vector space. We prove the main theorem regarding angles between collision subspaces in Section \ref{S3}. In Section \ref{S4} we apply this concept to billiards in the plane and make a comparison to existing billiard collision theorems. Section \ref{S5} provides a brief commentary on the limitations of $N$-body billiards and the techniques used in this paper, while also providing suggestions for further directions in the study of this problem.

\section{Billiards, collisions, and linear subspaces} \label{S2}
\subsection[N-body billiards]{$\boldsymbol{N}$-body billiards}

Motivated by the high-energy limit of the $N$-body problem, \cite{FKM} constructs \emph{N-body billiards}, a~formulation of an $N$-billiard dynamical system. We provide an outline of this system \mbox{below}. Consider an $N$ massive point particle system in $d$-dimensional Euclidean space $\mathbb{R}^{d}$ by its con\-fi\-gu\-ra\-tion space
\begin{gather*}
E = \big(\mathbb{R}^{d}\big)^N \approx \mathbb{R}^{N} \otimes \mathbb{R}^{d}.
\end{gather*}
Within $E$ there are $\binom{N}{2}$ \emph{binary collision subspaces}
\begin{gather*}
\Delta_{ij} = \big\{q=(q_1,\dots, q_N) \in \big(\mathbb{R}^{d}\big)^N\colon q_i = q_j \big\} \subset E
\end{gather*}
for some $i \neq j$. A \emph{billiard trajectory} will be a polygonal curve $\ell\colon \mathbb{R} \to E$, all of whose vertices are collisions (i.e., vertices of $\ell$ lie in $\Delta_{ij}$ for some $i$, $j$). When $\ell$ intersects a collision subspace $\Delta_{ij}$ it instantaneously changes direction by the law ``angle of incidence equals angle of reflection,'' given by equations \eqref{CE} and \eqref{CM} below. We call a \emph{collision point} a time $t^*$ for which $\ell(t^*) \in \Delta_{ij}$ for some distinct $1 \leq i, j \leq N$. We assume collision points are discrete and that no edge of $\ell$ lies within a collision subspace. The velocities $v_-$, $v_+$ of $\ell$ immediately before and after collision with~$\Delta_{ij}$ are locally constant and well-defined. These velocities undergo a jump $v_- \mapsto v_+$ at collision. Define
\begin{gather*}
\pi_{\Delta_{ij}}\colon \ E \to \Delta_{ij}
\end{gather*}
to be the orthogonal projection onto $\Delta_{ij}$. We require that each velocity jump follow the rules
\begin{gather}
\|v_-\|=\|v_+\|,
\label{CE}\\
\pi_{\Delta_{ij}}(v_-) = \pi_{\Delta_{ij}}(v_+),
\label{CM}
\end{gather}
which we consider as conservation of energy and conservation of linear momentum, respectively. Without loss of generality, we assume the billiard trajectory $\ell$ has unit speed.

An astute reader will notice that equation \eqref{CM} is ambiguous if the collision point $t^*$ is a time at which multiple collisions occur (e.g., triple collision or simultaneous binary collisions). This is analogous to trying to define standard billiard dynamics at a vertex of a polygonal billiard table. However, we only explore the case with $N=3$ point-billiards, and hence triple collision is the only scenario in which $t^*$ could be in more than one collision subspace. In~\cite{FKM} the issue is addressed by agreeing to choose one of the collision subspaces to which $t^*$ will belong and use only that subspace in implementing the conservation of momentum rule \eqref{CM}. In particular, billiard trajectories with multiple collisions include extra structure of labelling of collision points (see Section~2.1 of~\cite{FKM} for additional details).

This construction also leads to non-deterministic dynamics. For a given $v_- \in E \setminus \boldsymbol{0}$ and $t^* \in \Delta_{ij}$, if the binary collision subspaces are codimension $d$, $d\geq 1$, there is a $(d-1)$-dimensional sphere's worth of choices for outgoing velocities $v_+$. Even if $d=1$, the dynamics are non-deterministic, as the 0-sphere consists of two choices. It is standard to turn this case into a~deterministic process by requiring transversality: $v_+ \neq v_-$ at each collision. This is exactly the case for $N$ masses on a line. Similar degenerate billiard constructions have been studied~\cite{Bol2016,Bol2017} along with their connections to celestial mechanics.

Consider a simple nontrivial case, $N=1$ and $d=2$, and hence there are no binary collision subspaces. Suppose the point-billiard is launched into a planar wedge of angle measure $\alpha$. Can such a particle be trapped inside the wedge for infinite time? The answer is negative, and moreover through a simple geometric argument we can provide an upper bound on the total number of collisions of the particle with the sides of the wedge.

\begin{Theorem}
Consider a billiard trajectory inside a wedge with angle measure $\alpha$. The maximum number of collisions within the angle is $\big\lceil \frac{\pi}{\alpha}\big\rceil$ where $\lceil \cdot \rceil$ is the least integer function.
\end{Theorem}

\begin{figure}[ht]\centering
\includegraphics[scale=1]{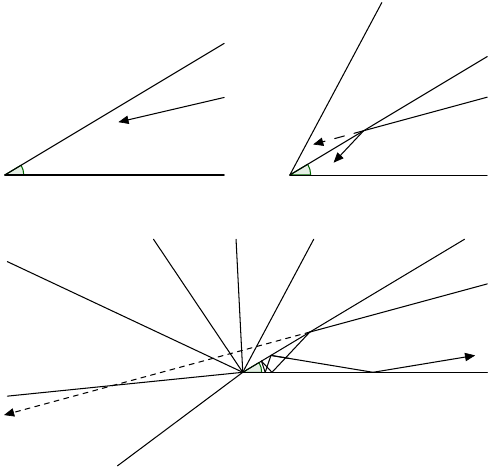}
\put(-180,128){\makebox(0,0)[lb]{$a$}}
\put(-53,128){\makebox(0,0)[lb]{$b$}}
\put(-120,-2){\makebox(0,0)[lb]{$c$}}
\caption{$(a)$ The billiard enters the wedge; $(b)$ upon colliding with one wall, reflect the wedge across the wall so the reflected trajectory is a straight line; $(c)$ repeat until the straight-line trajectory leaves the reflected wedge. The angle unfolds six times for a total of six collisions in the original wedge.}
\label{AngleUnfold}
\end{figure}

This can clearly be seen as follows. Consider an incoming billiard trajectory into the wedge of angle $\alpha$. Reflect the angle itself $\alpha$ across the side of impact and consider the image of the trajectory under those reflections. The billiard trajectory's image is a straight line, which can be seen to obey the rule ``angle of incidence equals angle of reflection''. This argument is illustrated in Fig.~\ref{AngleUnfold}.

Our aim is to use this technique in higher dimensions to bound the total number of collisions by using the collision subspaces as the ``walls'' of the wedge.

\subsection[Tensor construction of the N-billiard system]{Tensor construction of the $\boldsymbol N$-billiard system}

A useful tool in the $N$-body problem is the \emph{mass metric} on $\mathbb{R}^{N}$:
\begin{equation}
\langle \textbf{v}, \textbf{w} \rangle _M = \sum_{i=1}^N m_iv_iw_i
\label{massmetric}
\end{equation}
for vectors $\textbf{v},\textbf{w} \in \mathbb{R}^{N}$ and masses $m_i >0$. It follows from this definition that the squared norm is twice the kinetic energy.

In the configuration space $\mathbb{R}^{N} \otimes \mathbb{R}^{d}$, we use the mass metric on $\mathbb{R}^{N}$ and the standard Euclidean inner product on $\mathbb{R}^{d}$. Let $\varepsilon_i \in \mathbb{R}^{N}$ denote the $i^{\rm th}$ standard basis vector.
Define $E_i := \frac{\varepsilon_i}{\sqrt{m_i}}$ so that $E_i$ is a unit vector in $\mathbb{R}^{N}$ with respect to the mass metric.

It will be useful to translate the definition of the binary collision subspace into our tensor product construction as the following span of orthonormal elements of $\mathbb{R}^{N} \otimes \mathbb{R}^{d}$:
\begin{equation*}
 \Delta_{ij} = \left\{E_1 \otimes q_1, \dots, \frac{\varepsilon_i+\varepsilon_j}{\sqrt{m_i+m_j}} \otimes q_i, \dots ,E_N \otimes q_N\colon q_k \in \mathbb{R}^{d},\ \|q_k\|=1 \right\}.
\label{def:Delta:ij}
\end{equation*}

\begin{Remark}
In the collision subspaces, we adopt the convention that the ``location'' index will match that of the smaller of the two point mass indices, e.g., the element $\frac{\varepsilon_i + \varepsilon_j}{\sqrt{m_i+m_j}} \otimes q_i$ will be used instead of $\frac{\varepsilon_i + \varepsilon_j}{\sqrt{m_i+m_j}} \otimes q_j$ as a basis element in $\Delta_{ij}$. And we will continue to assume that each~$q_k$ is itself a unit vector in $\mathbb{R}^{d}$. We omit the word ``span'' and write subspaces $U = \{u_1,\dots, u_N\}$ to mean the $\mathbb{R}^{}$-linear span of the listed basis vectors. We shall also write these subspaces in terms of an orthonormal basis (even though an orthogonal basis is good enough).
\end{Remark}

\subsection{Linear algebra and principal angles}

Basic ideas from linear algebra provide a tool for computing the angle between linear subspaces of a vector space.

\begin{Definition}\label{padef1}
Let $F$, $G$ be subspaces of $\mathbb{R}^{n}$ with $\dim(F)=p$, $\dim(G)=q$ and let $r := \min\{p,q\}.$ The
\emph{principal angles} $\angle(F,G) = [\theta_1, \theta_2, \dots, \theta_r]$, are given by
\begin{gather*}
\cos(\theta_1) = \max_{\substack{u \in F \\ \|u\|=1}} \max_{\substack{v \in G \\ \|v\|=1}} \langle u, v\rangle = \langle u_1, v_1\rangle ,
\end{gather*}
where $\langle \cdot, \cdot\rangle $ is the standard Euclidean inner product. For each $1<k\leq r$,
\begin{gather*}
\cos(\theta_k) = \max_{\substack{u \in F \\ \|u\|=1 \\ \langle u, u_i\rangle =0}} \max_{\substack{v \in G \\ \|v\|=1 \\ \langle v, v_i\rangle =0}} \langle u, v\rangle = \langle u_k, v_k\rangle
\end{gather*}
for each $1 \leq i \leq k-1$. The vectors $u_k$, $v_k$ which realize the angle $\theta_k$ are called \emph{principal vectors}.
\end{Definition}

By construction, we have that $0 \leq \theta_1 \leq \dots \leq \theta_r \leq \pi/2$. An equivalent definition in terms of cosines of the principal angles can be found using the singular value decomposition of the matrix $\mathcal{F}^\top \mathcal{G}$ where $\mathcal{F}$, $\mathcal{G}$ are matrices whose columns are orthonormal bases for the subspaces~$F$ and $G$, respectively. See \cite{BG} for details.

The next lemma and example demonstrate that the angles between linear subspaces follow similar properties to what one would expect in the standard Euclidean geometry.

\begin{Lemma}[\cite{AK2007,AJK2010}]
Let $F$, $G$ be subspaces of $\mathbb{R}^{n}$ with $\dim(F) =p$, $\dim(G) =q$. Furthermore let $\angle(F,G)$ and $\angle^*(F,G)$ denote the vector of principal angles between $F$ and $G$ in increasing order and decreasing order, respectively.
\begin{enumerate}\itemsep=0pt
\item[$1.$] $\angle(F,G) = \angle(G,F)$.
\item[$2.$] $[0,\dots, 0, \angle(F,G)] = \big[0,\dots,0, \angle\big(F^\perp,G^\perp\big)\big]$, with $\max\{n-p-q,0\}$ zeros on the left and $\max\{p+q-n,0\}$ zeros on the right.
\item[$3.$] $\big[0,\dots,0, \angle\big(F,G^\perp\big)\big] = \big[0,\dots,0,\angle\big(F^\perp,G\big)\big]$, with $\max\{q-p,0\}$ zeros on the left and $\max\{p-q,0\}$ zeros on the right.
\item[$4.$] $\big[\angle(F,G), \frac{\pi}{2},\dots,\frac{\pi}{2}\big] = \big[0,\dots,0, \frac{\pi}{2}-\angle^*\big(F,G^\perp\big)\big]$, with $\max\{p-q,0\}$ $\frac{\pi}{2}$'s on the left and $\max\{p+q-n,0\}$ zeros on the right.
\end{enumerate}
\end{Lemma}

The first property follows from Definition~\ref{padef1}, and we do not provide proofs for the rest of~the properties. Proofs can be found as Theorem~2.7 in \cite{AK2007} or as Theorems 2.6 and 2.7 in \cite{AJK2010}.

\begin{Example}
Consider subspaces $L$ and $M$ in $\mathbb{R}^{6}$.
\begin{itemize}\itemsep=0pt
\item Suppose $\dim(L)=\dim(M)=2$ and $\angle(L,M) = \big[\frac{\pi}{3},\frac{\pi}{2}\big]$. Then $\angle\big(L^\perp,M^\perp\big) = \big[0,0,\frac{\pi}{3},\frac{\pi}{2}\big]$, and $\angle\big(L,M^\perp\big) = \angle\big(L^\perp,M\big) = \big[0,\frac{\pi}{6}\big]$.

\item Suppose $\dim(L)=\dim(M)=4$ and $\angle(L,M) = \big[0,0,\frac{\pi}{4},\frac{\pi}{3}\big].$ We conclude $\angle\big(L^\perp,M^\perp\big) = \big[\frac{\pi}{4},\frac{\pi}{3}\big]$ and $\angle\big(L,M^\perp\big) = \angle\big(L,M^\perp\big) = \big[\frac{\pi}{6},\frac{\pi}{4}\big]$.
\end{itemize}
\end{Example}

Though a subtlety that isn't obvious in the two examples above is that, given $\angle(L,M)$, the~angles that appear in the vector $\angle\big(L^\perp,M^\perp\big)$ are taken from the list $\angle(L,M)$ from largest to smallest. We illustrate this point in the next example.

\begin{Example}
If $\dim(L)=\dim(M) = 8$ in $\mathbb{R}^{13}$ and
\begin{gather*}
\angle(L,M) = \bigg[0,\frac{\pi}{6},\frac{\pi}{6},\frac{\pi}{4},
\frac{\pi}{3},\frac{\pi}{3},\frac{\pi}{3},\frac{\pi}{2}\bigg],
\end{gather*}
then
\begin{gather*}
\angle(L^\perp,M^\perp) = \bigg[\frac{\pi}{4},\frac{\pi}{3},\frac{\pi}{3},\frac{\pi}{3},\frac{\pi}{2}\bigg],
\end{gather*}
which are the five largest angles in the vector $\angle(L,M)$.
\end{Example}

\begin{Corollary}If $U,V \subset \mathbb{R}^{n}$ are codimension 1 subspaces, then the nonzero angle between $U$ and $V$ is $\angle\big(U^\perp, V^\perp\big)$. That is, the nonzero angle between these subspaces is precisely the angle between their normal vectors.
\end{Corollary}

Our first goal is to compute the angles between the collision subspaces $\Delta_{ij}$ and $\Delta_{kl}$ for some positive integers $i$, $j$, $k$, $l$ and $1 \leq i,j,k,l \leq N$. By our definitions of principal angles, there will be $d(N-1)$ principal angles between the codimension $d$ collision subspaces.

\section{Angles between collision subspaces}\label{S3}

Suppose the masses in each of the $N$ point billiards are $m_i$, $m_j$, $m_k$, $m_l$ for positive integers $i$, $j$, $k$, $l$. We can compute the principal angles between the binary collision subspaces in terms of~their respective masses.

\begin{Theorem}\label{thm2}
Let $i$, $j$, $k$, $l$ be distinct integers satisfying $1 \leq i,j,k,l \leq N$.
\begin{enumerate}\itemsep=0pt
\item[$a)$] The first $d(N-2)$ principal angles between $\Delta_{ij}$ and $\Delta_{kl}$ are $0$ and the remaining $d$ principal angles are $\pi/2$.
\item[$b)$] The first $d(N-2)$ principal angles between $\Delta_{ij}$ and $\Delta_{jk}$ are $0$ and the last $d$ principal angles are
\begin{equation*}
\theta = \arccos\bigg(\dfrac{2(m_i+m_k)-m_j}{\sqrt{(m_i+m_j+4m_k)(4m_i+m_j+m_k)}} \bigg).
\label{angleeqn}
\end{equation*}
\end{enumerate}
\end{Theorem}

\begin{proof}First, by Definition \ref{padef1} we note that the number of principal angles which are 0 bet\-ween~$\Delta_{ij}$ and $\Delta_{kl}$ is exactly the dimension of their intersection. In both cases ($a$) and ($b$), the dimen\-sion of the intersection is $d(N-2)$.

\begin{figure}[hpbt]\centering
\includegraphics[scale=1]{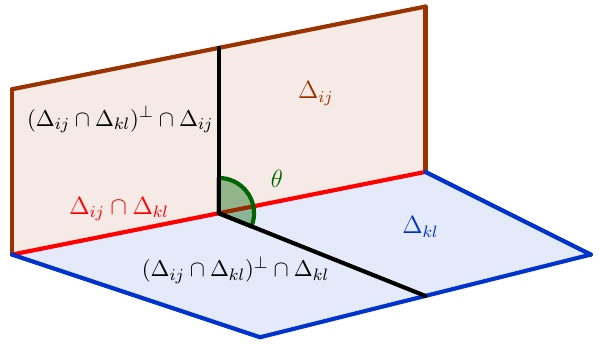}
\caption{The geometric angle $\theta$ computed between collision subspaces.}\label{AngleBetSubspaces}
\end{figure}

Writing the collision subspaces as the span of basis elements in tensor product form, we see that
\begin{gather*}
\Delta_{ij} = \bigg\{E_1 \otimes q_1, \dots, \frac{\varepsilon_i+\varepsilon_j}{\sqrt{m_i+m_j}} \otimes q_i, \dots, E_N \otimes q_N \bigg\}
\end{gather*}
and
\begin{gather*}
\Delta_{kl} = \bigg\{E_1 \otimes q_1, \dots, \frac{\varepsilon_k+\varepsilon_l}{\sqrt{m_k+m_l}} \otimes q_k, \dots, E_N \otimes q_N \bigg\}.
\end{gather*}

Geometrically the desired angle between the subspaces can be found by computing the angle between $(\Delta_{ij} \cap \Delta_{kl})^\perp \cap \Delta_{ij}$ and $(\Delta_{ij}\cap \Delta_{kl})^\perp \cap \Delta_{kl}$. See Fig.~\ref{AngleBetSubspaces}.
First, we find that
\begin{gather*}
\Delta_{ij}\cap \Delta_{kl} = \bigg\{ E_1\otimes q_1, \dots, \frac{\varepsilon_i + \varepsilon_j}{\sqrt{m_i+m_j}} \otimes q_i, \dots, \frac{\varepsilon_k + \varepsilon_l}{\sqrt{m_k+m_l}} \otimes q_k, \dots, E_N \otimes q_N \bigg\},
\end{gather*}
and hence
\begin{gather*}(\Delta_{ij}\cap\Delta_{kl})^\perp = \bigg\{\frac{m_j\varepsilon_i -m_i\varepsilon_j}{\sqrt{m_im_j(m_i+m_j)}} \otimes \beta_1, \frac{m_l\varepsilon_k-m_k\varepsilon_l}{\sqrt{m_km_l(m_k+m_l)}}\otimes \beta_2 \bigg\}
\end{gather*}
for some arbitrary unit vectors $\beta_1, \beta_2 \in \mathbb{R}^{d}$. Then
\begin{gather*}
(\Delta_{ij}\cap \Delta_{kl})^\perp \cap \Delta_{ij} = \bigg\{ \dfrac{m_l\varepsilon_k-m_k\varepsilon_l}{\sqrt{m_km_l(m_k+m_l)}} \otimes \beta_2 \bigg\}
\end{gather*}
and
\begin{gather*}(\Delta_{ij}\cap \Delta_{kl})^\perp \cap \Delta_{kl} = \bigg\{ \dfrac{m_j\varepsilon_i-m_i\varepsilon_j}{\sqrt{m_im_j(m_i+m_j)}} \otimes \beta_1 \bigg\}.
\end{gather*}
Thus
\begin{gather*}
\big[(\Delta_{ij}\cap \Delta_{kl})^\perp \cap \Delta_{ij}\big] \perp \big[ (\Delta_{ij}\cap \Delta_{kl})^\perp \cap \Delta_{kl} \big].
\end{gather*}
Therefore regardless of the individual masses, $\Delta_{ij}$ and $\Delta_{kl}$ are orthogonal to one another when $i$, $j$, $k$, $l$ are distinct.

We now turn our attention to the case where one of the indices is the same across the two collision subspaces. Repeating the same calculations as above, we find that
\begin{gather*}
\Delta_{ij}\cap \Delta_{jk} = \bigg\{ E_1\otimes q_1, \dots, \dfrac{\varepsilon_i + \varepsilon_j + \varepsilon_k}{\sqrt{m_i+m_j+m_k}} \otimes q_i, \dots, E_N \otimes q_N\bigg\}
\end{gather*}
and
\begin{gather*}
(\Delta_{ij}\cap \Delta_{jk})^\perp = \bigg\{ \frac{m_j\varepsilon_i-m_i\varepsilon_j}{\sqrt{m_i+m_j}} \otimes \beta_1, \dfrac{\varepsilon_i+\varepsilon_j - 2\varepsilon_k}{\sqrt{m_i+m_j+4m_k}} \otimes \beta_3\bigg\}
\end{gather*}
for an arbitrary unit vector $\beta_3 \in \mathbb{R}^{d}$.
Then
\begin{gather*}
(\Delta_{ij}\cap \Delta_{jk})^\perp \cap \Delta_{ij} = \bigg\{ \dfrac{\varepsilon_i+\varepsilon_j - 2\varepsilon_k}{\sqrt{m_i+m_j+4m_k}} \otimes \beta_3\bigg\}
\end{gather*}
and
\begin{gather*}
(\Delta_{ij}\cap \Delta_{jk})^\perp \cap \Delta_{jk} = \bigg\{ \dfrac{2\varepsilon_i-\varepsilon_j - \varepsilon_k}{\sqrt{4m_i+m_j+m_k}} \otimes \beta_3\bigg\}.
\end{gather*}

Therefore by Definition \ref{padef1}, the nonzero angle between these two subspaces is given by
\begin{equation*}
	\cos(\theta) = \dfrac{2(m_i+m_k)-m_j}{\sqrt{(m_i+m_j+4m_k)(4m_i+m_j+m_k)}}.\tag*{\qed}
 \end{equation*}
 \renewcommand{\qed}{}
\end{proof}

If the masses are all equal, the result is easy to state.

\begin{Corollary}
Let $i$, $j$, $k$, $l$ be distinct integers satisfying $1 \leq i,j,k,l \leq N$. Suppose $m_i = m_j = m_k = m_l>0$.
\begin{enumerate}\itemsep=0pt
\item[$a)$] The first $d(N-2)$ principal angles between $\Delta_{ij}$ and $\Delta_{kl}$ are $0$ and the remaining $d$ principal angles are $\pi/2$.
\item[$b)$] The first $d(N-2)$ principal angles between $\Delta_{ij}$ and $\Delta_{jk}$ are $0$ and the remaining $d$ principal angles are $\pi/3$.
\end{enumerate}
\end{Corollary}

\section{Billiard trajectories and collision bounds in the plane}\label{S4}

\subsection{A primer on Jacobi coordinates and the mass metric}

We follow the approach of Sections 3 and 7 of \cite{M1}. We consider the planar 3-billiard ball problem whose configuration space is $\mathbb{C}^3$. A vector $\textbf{q} = (q_1,q_2,q_3) \in \mathbb{C}^3$ represents a located triangle with each of its components representing the vertices of the triangle.

The \emph{mass metric} on the configuration space $\mathbb{C}^3$ is the Hermitian inner product
\begin{gather*}
\langle \textbf{v}, \textbf{w}\rangle _M = m_1 \overline{v_1}w_1 + m_2 \overline{v_2}w_2 +m_3 \overline{v_3}w_3.
\end{gather*}
This is consistent with equation \eqref{massmetric} used in the tensor construction.

A translation of this located triangle \textbf{q} by $c \in \mathbb{C}$ is given by the located triangle $\textbf{q}+c\textbf{1}$ where $\textbf{1}= (1,1,1)$. Define
\begin{gather*}
\mathbb{C}_0^2:=\textbf{1}^\perp = \big\{\textbf{q}\in\mathbb{C}^3\colon m_1q_1+m_2q_2+m_3q_3 =0\big\}
\end{gather*}
to be the set of planar three-body configurations whose center of mass $\textbf{q}_{cm}$ is at the origin. This two-dimensional complex space represents the quotient space of $\mathbb{C}^3$ by translations.

\begin{Definition}
The \emph{Jacobi coordinates} for $\mathbb{C}_0^2 :=\big\{\textbf{q} \in \mathbb{C}^3\colon \textbf{q}_{cm}=0\big\}$ are given by
\begin{gather*}
v = \mu_1(q_1-q_2), \qquad w = \mu_2\bigg(q_3 - \frac{m_1q_1+m_2q_2}{m_1+m_2}\bigg),
\end{gather*}
where $\frac{1}{\mu_1^2} = \frac{1}{m_1}+\frac{1}{m_2}$ and $\frac{1}{\mu_2^2}= \frac{1}{m_3}+\frac{1}{m_1+m_2}$.
\end{Definition}

These are normalized coordinates which diagonalize the restriction of the mass metric to $\mathbb{C}_0^2$, see Fig.~\ref{Fig3}(a). From this we can define the complex linear projection
\begin{gather*}
\pi_{\rm tr}\colon\ \mathbb{C}^3 \to \mathbb{C}^2, \qquad (q_1,q_2,q_3) \mapsto (v,w),
\end{gather*}
which realizes the metric quotient of $\mathbb{C}^3$ by translations.

It is worthwhile to note that using Jacobi coordinates and our map $\pi_{\rm tr}$, all of the triple collision triangles $(q,q,q) \in \mathbb{C}^3$ are mapped to the origin.

\subsection{Collision bounds in the plane}

Consider equal masses $M$ in the plane. This changes the linear projection into
\begin{gather*}
\pi_{\rm tr}\colon\ \mathbb{C}^3 \to \mathbb{C}^2 , \qquad
(q_1,q_2,q_3) \mapsto \sqrt{M}\bigg(\frac{1}{\sqrt{2}}(q_1-q_2), \sqrt{\frac{2}{3}}\bigg(q_3-\frac{1}{2}(q_1+q_2) \bigg) \bigg).
\end{gather*}
The mass $M$ acts as a dilation factor, so we assume the mass to be unit henceforth.

The codimension 2 binary collision subspaces can be defined in terms of the complex coordinates as follows:
\begin{gather*}
\Delta_{12} = \big\{(q_1,q_2,q_3)\in \mathbb{C}^3 \colon q_1 = q_2\big\},
\\
\Delta_{23} = \big\{(q_1,q_2,q_3)\in \mathbb{C}^3 \colon q_2 = q_3\big\},
\\
\Delta_{13} = \big\{(q_1,q_2,q_3)\in \mathbb{C}^3 \colon q_1 = q_3\big\}.
\end{gather*}

\begin{Example} Using Theorem \ref{thm2}, we can see that $\angle(\Delta_{12}, \Delta_{23}) = \big[0,0,\frac{\pi}{3},\frac{\pi}{3}\big]$. In this case the principal vectors are
\begin{gather*}
\theta_1 =0\colon\ \frac{1}{\sqrt{3}}(1,1,1)\qquad \text{and itself},
\\
\theta_2 =0\colon\ \frac{1}{\sqrt{3}}(i,i,i)\qquad \text{and itself},
\\
\theta_3 = \frac{\pi}{3}\colon\ \frac{1}{\sqrt{6}}(-1,-1,2)\qquad \text{and}\qquad \frac{1}{\sqrt{6}}(-2,1,1),
\\
\theta_4 = \frac{\pi}{3}\colon\ \frac{1}{\sqrt{6}}(-i,-i,2i)\qquad \text{and} \qquad \frac{1}{\sqrt{6}}(-2i,i,i).
\end{gather*}
\end{Example}

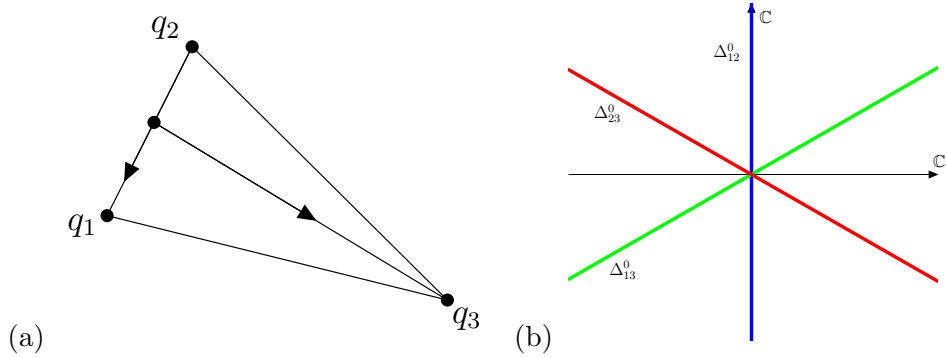
\begin{figure}[htbp]
\centering
\includegraphics[scale=1]{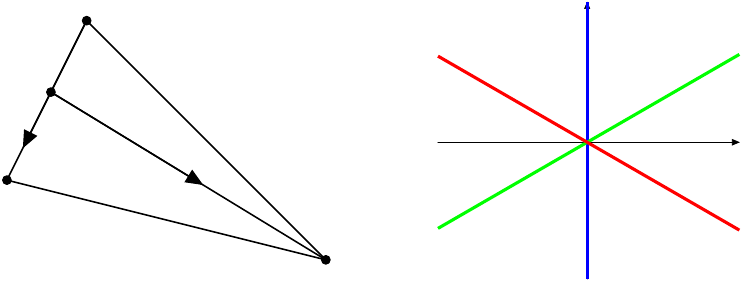}
\put(-73,128){\makebox(0,0)[lb]{$\mathbb{C}$}}
\put(-0,65){\makebox(0,0)[lb]{$\mathbb{C}$}}
\put(-98,110){\makebox(0,0)[lb]{$\Delta_{12}^0$}}
\put(-136,80){\makebox(0,0)[lb]{$\Delta_{23}^0$}}
\put(-122,30){\makebox(0,0)[lb]{$\Delta_{23}^0$}}
\put(-197,7){\makebox(0,0)[lb]{$q_3$}}
\put(-329,130){\makebox(0,0)[lb]{$q_2$}}
\put(-368,49){\makebox(0,0)[lb]{$q_1$}}
\put(-290,-11){\makebox(0,0)[lb]{$a$}}
\put(-80,-11){\makebox(0,0)[lb]{$b$}}
\caption{$(a)$ Jacobi vectors; $(b)$ an ``overhead'' view of $\mathbb{C}_0^2$ and the collision subspaces.}\label{Fig3}
\end{figure}

The image of these subspaces under $\pi_{\rm tr}$ are
\begin{gather*}
\Delta_{12}^0 = \big\{(v,w)\in \mathbb{C}_0^2 \colon v=0\big\} = \operatorname{span}_\mathbb{C}\{(0,1)\},
\\
\Delta_{23}^0= \bigg\{(v,w)\in \mathbb{C}_0^2 \colon w=-\frac{1}{\sqrt{3}}v\bigg\} = \operatorname{span}_\mathbb{C} \bigg\{\bigg(\frac{\sqrt{3}}{2},-\frac{1}{2}\bigg)\bigg\},
\\
\Delta_{13}^0 = \bigg\{(v,w)\in \mathbb{C}_0^2 \colon w=\frac{1}{\sqrt{3}}v\bigg\} = \operatorname{span}_\mathbb{C} \bigg\{\bigg(\frac{\sqrt{3}}{2},\frac{1}{2}\bigg)\bigg\}.
\end{gather*}
Each of these codimension 2 subspaces are planes in $\mathbb{C}_0^2$, as pictured in Fig.~\ref{Fig3}(b).

Through this reduction via Jacobi coordinates, the angles between these subspaces are $\angle\big(\Delta_{12}^0,\Delta_{23}^0\big) = \angle\big(\Delta_{12}^0,\Delta_{13}^0\big) = \angle\big(\Delta_{23}^0,\Delta_{13}^0\big) =\big[\frac{\pi}{3}, \frac{\pi}{3}\big].$

And to further our previous example, we can observe that the principal vectors for
$\angle\big(\Delta_{12}^0{,} \Delta_{23}^0\big)\!$ are $(0,1)$ and $\big({-}\frac{\sqrt{3}}{2}, \frac{1}{2}\big)$ for
the first angle, and $(0,i)$ and $\big({-}\frac{\sqrt{3}}{2}i, \frac{1}{2}i\big)$ for the second angle. But in
fact one can check that indeed
\begin{gather*}
\pi_{\rm tr}\bigg(\bigg({-}\frac{1}{\sqrt{6}}, -\frac{1}{\sqrt{6}}, \frac{2}{\sqrt{6}}\bigg)\bigg) = (0,1)
\end{gather*}
and
\begin{gather*}
\pi_{\rm tr}\bigg(\bigg({-}\frac{2}{\sqrt{6}}, \frac{1}{\sqrt{6}}, \frac{1}{\sqrt{6}}\bigg)\bigg) = \bigg({-}\frac{\sqrt{3}}{2},\frac{1}{2}\bigg).
\end{gather*}

Since $\pi_{\rm tr}$ is linear, we know that the image of the other pair of principal vectors between $\Delta_{12}$ and $\Delta_{23}$ will also be the principal vectors between $\Delta_{12}^0$ and $\Delta_{23}^0$. That is, the image of a nonzero principal vector under $\pi_{\rm tr}$ is still a principal vector.

\begin{Theorem}\label{colthm1}
In the equal mass planar $3$-body billiard problem there can be at most $3$ collisions.
\end{Theorem}

To prove the theorem, we need to following lemma.

\begin{Lemma}
Let $V_1$ and $V_2$ be arbitrary vectors in two collision subspaces in $\mathbb{C}_0^2$ and let $\hat{\theta}$ denote the angle between the vectors $V_1$ and $V_2$. Then $\frac{\pi}{3} \leq \hat{\theta} \leq \frac{\pi}{2}$.
\end{Lemma}

\begin{proof} Without loss of generality we consider two of our three collision subspaces, namely $\Delta_{12}^0$ and $\Delta_{23}^0$ and consider two arbitrary vectors $V_1$ and $V_2$ in $\Delta_{12}^0$ and $\Delta_{23}^0$, respectively (see Fig.~\ref{Fig4}).


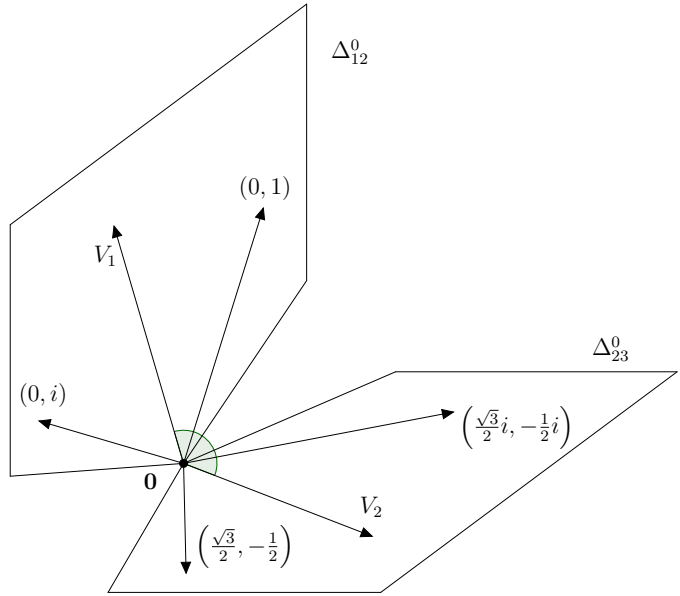
\begin{figure}[hbtp]\centering
\includegraphics[scale=.97]{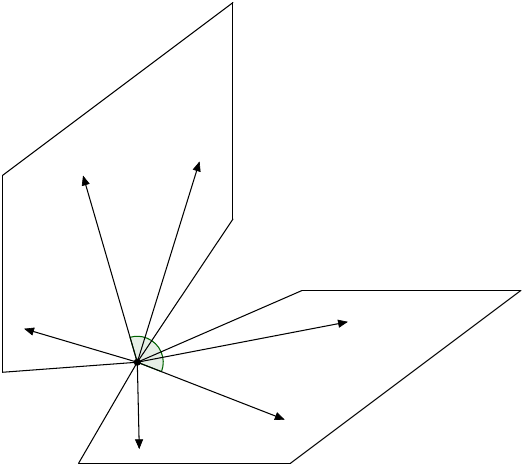}
\put(-130,190){\makebox(0,0)[lb]{$\Delta_{12}^0$}}
\put(-30,87){\makebox(0,0)[lb]{$\Delta_{12}^0$}}
\put(-220,125){\makebox(0,0)[lb]{$V_1$}}
\put(-241,68){\makebox(0,0)[lb]{$(0,i)$}}
\put(-163,145){\makebox(0,0)[lb]{$(0,1)$}}
\put(-115,25){\makebox(0,0)[lb]{$V_2$}}
\put(-82,60){\makebox(0,0)[lb]{$\left(\frac{\sqrt3}{2}i,-\frac12i\right)$}}
\put(-178,8){\makebox(0,0)[lb]{$\left(\frac{\sqrt3}{2},-\frac12\right)$}}
\put(-193,38){\makebox(0,0)[lb]{$\mathbf{0}$}}
\caption{The two subspaces $\Delta_{12}^0$ and $\Delta_{23}^0$ along with their orthonormal basis vectors. In our proof we look to measure the angle between arbitrary vectors $V_1 \in \Delta_{12}^0$ and $V_2 \in \Delta_{23}^0$.}\label{Fig4}
\end{figure}

Recall the definition of principal angles:
\begin{gather*}
\cos(\theta_1) = \max_{\substack{u \in \Delta_{12}^0 \\ \|u\|=1}} \max_{\substack{v \in \Delta_{23}^0 \\ \|v\|=1}} \langle u, v\rangle := \langle u_1, v_1\rangle ,
\end{gather*}
where $u_1$ and $v_1$ are the principal vectors which realize this principal angle.

From our earlier calculations, we know that $\cos(\theta_1) = \frac{\pi}{3}$ and that our principal angles always satisfy $0 \leq \theta_i \leq \frac{\pi}{2}$. Because $\cos(\theta)$ is a decreasing function on the interval $0 \leq \theta \leq \frac{\pi}{2}$, we see that
\begin{gather*}
\cos(\theta_1) = \langle u_1, v_1\rangle \geq \langle V_1, V_2\rangle = \cos\big(\hat{\theta}\big),
\end{gather*}
because the inner product is maximized. Hence the possible angles between the vectors $V_1$ and~$V_2$ must satisfy $\frac{\pi}{3} \leq \hat{\theta} \leq \frac{\pi}{2}.$ The argument and calculation is the same if we choose any pair of these collision subspaces. This proves the lemma.

In fact, the angle $\hat{\theta} = \frac{\pi}{2}$ can be realized if we let $V_1 = \big(0,\frac{\sqrt{3}}{2}-\frac{1}{2}i\big)$ and $V_2 = \big(\frac{\sqrt{3}+3i}{4},\frac{1+\sqrt{3}i}{4}\big)$.
\end{proof}

Using the preceding lemma, the proof of the theorem is short and follows from the ``unfolding the angle'' argument. We also refer to a sequence of collisions by the order in which the binary collisions occur. For example, the collision sequence (12)(23) indicates that $\ell(t)$ has a collision point in $\Delta_{12}$ first and then $\Delta_{23}$ second.

\begin{proof}
Consider an arbitrary piecewise linear trajectory $\ell(t)$ in $\mathbb{C}_0^2$, we aim to maximize the number of collisions. Without loss of generality, assume $\ell(t)$ intersects $\Delta_{12}^0$ first.
Suppose~$\ell(t)$ intersects $\Delta_{12}^0$ at time $t=t_1$ and let $V_1$ be a vector in $\Delta_{12}^0$ whose endpoint is this point of~intersection, $\ell(t_1)$. The next subspace $\ell(t)$ intersects can be either $\Delta_{23}^0$ or $\Delta_{13}^0$. Assume~$\ell(t)$ next visits $\Delta_{23}^0$ at time $t=t_2$, and let the vector $V_2 \in \Delta_{23}^0$ be a vector whose endpoint is this second point of intersection, $\ell(t_2)$. We know from the preceding lemma that the angle $\hat{\theta_1}$ between $V_1$ and $V_2$ satisfies $\frac{\pi}{3} \leq \hat{\theta_1}\leq \frac{\pi}{2}$.

\begin{figure}[hbtp]\centering
\includegraphics[scale=.95]{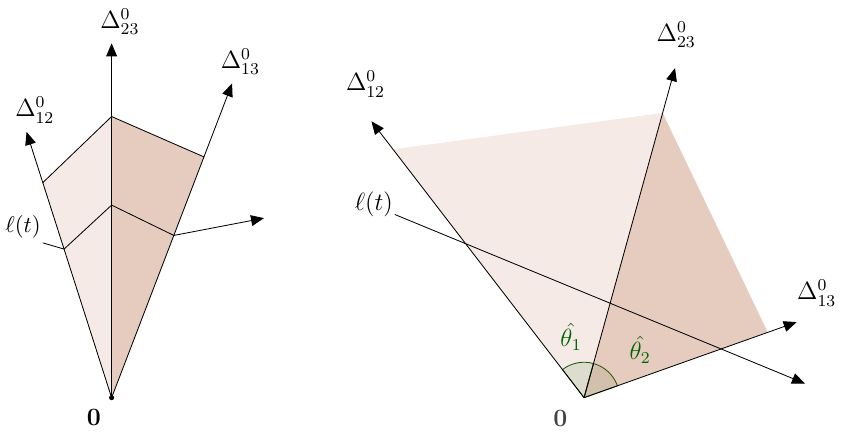}
\caption{The unfolding of the (12)(23)(13) trajectory.}\label{Fig5}
\end{figure}

We now repeat this process again. From $\ell(t_2)$, the trajectory $\ell(t)$ can now travel to either~$\Delta_{23}^0$ or $\Delta_{13}^0$. Without loss of generality, assume $\Delta_{13}^0$ is the next subspace. Let $\ell(t)$ intersect $\Delta_{13}^0$ at $\ell(t_3)$ for some time $t=t_3$ and let $V_3$ be a vector in $\Delta_{13}^0$ whose endpoint is at the point of intersection $\ell(t_3)$. Applying our lemma again, the angle $\hat{\theta_2}$ between $V_2$ and $V_3$ satisfies $\frac{\pi}{3} \leq \hat{\theta_2} \leq \frac{\pi}{2}$.
However, after time $t_3$, $\ell(t)$ cannot intersect any more collision subspaces. At best, $\hat{\theta_1}+\hat{\theta_2} = \frac{2\pi}{3}$, and any third angle will add at least $\frac{\pi}{3}$ by the previous lemma. So if we glue together the sectors spanned by $V_1$ and $V_2$, and $V_2$ and $V_3$ and flatten this angle, by ``unfolding the angle'' the trajectory $\ell(t)$ can intersect no more subspaces (see Fig.~\ref{Fig5}). This leaves us with a collision bound of at most $\big\lceil \frac{\pi}{(\pi/3)}\big\rceil =3$ possible collisions.
\end{proof}

\subsection{An arbitrary mass collision bound theorem}

Considering arbitrary masses and reusing the proof of Theorem \ref{colthm1}, we state the following:

\begin{Theorem}
For three arbitrary point-masses $m_i$, $m_j$, $m_k$, the maximum number of colli\-si\-ons~is
\begin{gather*}
\bigg\lceil \pi /\arccos\bigg( \dfrac{2(m_i+m_k)-m_j}{\sqrt{(m_i+m_j+4m_k)(4m_i+m_j+m_k)}} \bigg) \bigg\rceil.
\end{gather*}
\end{Theorem}

This expression can be simplified slightly. The expression above is symmetric in $m_i$ and $m_k$, so write $m_i = \alpha m_j$ and $m_k = \beta m_j$. The expression then is not directly dependent upon the masses but on the relative ratios of the masses, $\alpha$ and $\beta$:
\begin{gather*}
\bigg\lceil \pi / \arccos\bigg( \dfrac{2(\alpha+\beta)-1}{\sqrt{(\alpha+1 +4\beta)(4\alpha+1+\beta)}} \bigg) \bigg \rceil.
\end{gather*}

\begin{figure}[phtb]\centering
\includegraphics[scale=.9]{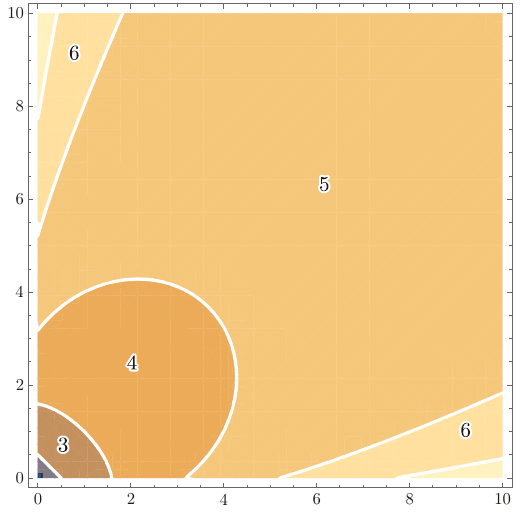}
\put(-130,-7){\makebox(0,0)[lb]{$\alpha$}}
\put(-260,122){\makebox(0,0)[lb]{$\beta$}}
\caption{A contour plot for the maximum number of collisions when $(\alpha, \beta) \in (0,10]\times(0,10]$ and $m_j =1$.}
\label{colplot}
\end{figure}

This expression also provides an interesting bound on the number of collisions. As seen in Fig.~\ref{colplot}, the number of collisions only seems to change when $\beta$ is large with $\alpha \ll \beta$ or vice-versa.


\subsection{A detour: the Foch sequence}

\begin{figure}[phtb]\centering
\includegraphics[scale=.9]{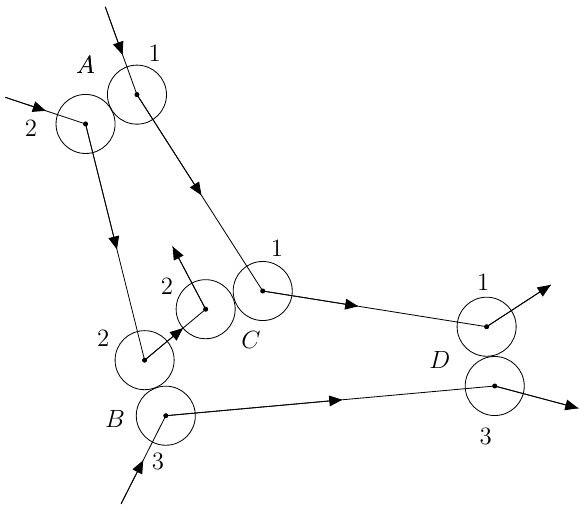}
\caption{The Foch Sequence (12)(23)(12)(13). After the (23) collision at $B$, particle 3 has to ``maneuver around'' particle 2 to collide with particle 1 in the fourth collision. See \cite{BD19} for an alternate viewpoint on this sequence. }
\label{FochSequence}
\end{figure}

A result by Murphy and Cohen \cite{MC1,MC2} states that the bound of the 3-billiard ball (seen as spheres) problem in $\mathbb{R}^{d}$ is four. The four-collision sequence (12)(23)(12)(13) is called the \emph{Foch sequence}, see Fig.~\ref{FochSequence}. Their proof is geometric and makes conditions on the locations of one of the balls in terms of the radii of the other billiard balls. However, it is interesting to note that despite how it may look at an initial glance, this Foch sequence does not contradict Theorem~\ref{colthm1}. In $N$-body billiards, the billiard balls are seen as point masses, and so no such considerations are necessary. In fact, if the radii shrink to zero and follow the details of the Murphy and Cohen proof, the Foch sequence is no longer possible, and the collision bound jumps from 4 to 3 when the radii reach zero.

\section{Future work and next steps}\label{S5}
The work in this paper is the result of attempts to solve the original problem: bounding the number of collisions in an $N$-billiard system using the computed angles between collision subspaces. The main collision theorems were possible due to the symplectic reduction using Jacobi coordinates, which reduced the angles between the reduced collision subspaces to all be nonzero. When $N>3$ and $d>2$, the above reduction techniques do not completely eliminate all nonzero angles. Further reductions should be considered, e.g., reducing by rotations.
This model is also inherently simplistic. Geometric and physical considerations are easily ignored (e.g., when considering hard rods on a line which cannot pass one another), but the bounds may still \mbox{exist}. Is~there another interpretation of the system that more closely matches a physical system? A~more rigorous study of this model to include such constraints is a logical next step.

\subsection*{Acknowledgements}
The author is grateful for the detailed comments and suggestions of the referees. This material is based upon work supported by the National Science Foundation under Grant \# DMS-1440140 while the author was in residence at the Mathematical Sciences Research Institute in Berkeley, California, during the Fall 2018 semester. This work was also supported by Discovery Project \# DP190101838, \emph{Billiards within confocal quadrics and beyond} from the Australian Research Council.

\pdfbookmark[1]{References}{ref}
\LastPageEnding


\begin{thebibliography}{99}
\footnotesize\itemsep=0pt

\bibitem{AT2018}
Albers P., Tabachnikov S., Introducing symplectic billiards, \href{https://doi.org/10.1016/j.aim.2018.05.037}{\textit{Adv.
 Math.}} \textbf{333} (2018), 822--867, \href{https://arxiv.org/abs/1708.07395}{arXiv:1708.07395}.

\bibitem{BMT2020}
Bialy M., Mironov A.E., Tabachnikov S., Wire billiards, the first steps,
 \href{https://doi.org/10.1016/j.aim.2020.107154}{\textit{Adv. Math.}} \textbf{368} (2020), 107154, 27~pages,
 \href{https://arxiv.org/abs/1905.13617}{arXiv:1905.13617}.

\bibitem{BG}
Bj\"orck A., Golub G.H., Numerical methods for computing angles between linear
 subspaces, \href{https://doi.org/10.2307/2005662}{\textit{Math. Comp.}} \textbf{27} (1973), 579--594.

\bibitem{Bol2016}
Bolotin S.V., Degenerate billiards, \href{https://doi.org/10.1134/S0081543816080046}{\textit{Proc. Steklov Inst. Math.}}
 \textbf{295} (2016), 45--62, \href{https://arxiv.org/abs/1606.06708}{arXiv:1606.06708}.

\bibitem{Bol2017}
Bolotin S.V., Degenerate billiards in celestial mechanics, \href{https://doi.org/10.1134/S1560354717010038}{\textit{Regul.
 Chaotic Dyn.}} \textbf{22} (2017), 27--53, \href{https://arxiv.org/abs/1612.08907}{arXiv:1612.08907}.

\bibitem{BFK1}
Burago D., Ferleger S., Kononenko A., Uniform estimates on the number of
 collisions in semi-dispersing billiards, \href{https://doi.org/10.2307/120962}{\textit{Ann. of Math.}} \textbf{147}
 (1998), 695--708.

\bibitem{BFK2}
Burago D., Ferleger S., Kononenko A., A geometric approach to semi-dispersing
 billiards, in Hard Ball Systems and the {L}orentz Gas, \href{https://doi.org/10.1007/978-3-662-04062-1_2}{\textit{Encyclopaedia
 Math. Sci.}}, Vol.~101, Springer, Berlin, 2000, 9--27.

\bibitem{BD19}
Burdzy K., Duarte M., A lower bound for the number of elastic collisions,
 \href{https://doi.org/10.1007/s00220-019-03399-3}{\textit{Comm. Math. Phys.}} \textbf{372} (2019), 679--711,
 \href{https://arxiv.org/abs/1803.00979}{arXiv:1803.00979}.

\bibitem{FKM}
F\'ejoz J., Knauf A., Montgomery R., Lagrangian relations and linear point
 billiards, \href{https://doi.org/10.1088/1361-6544/aa5b26}{\textit{Nonlinearity}} \textbf{30} (2017), 1326--1355,
 \href{https://arxiv.org/abs/1606.01420}{arXiv:1606.01420}.

\bibitem{Gal2003}
Galperin G., Playing pool with {$\pi$} (the number {$\pi$} from a billiard
 point of view), \href{https://doi.org/10.1070/RD2003v008n04ABEH000252}{\textit{Regul. Chaotic Dyn.}} \textbf{8} (2003), 375--394.

\bibitem{AK2007}
Knyazev A., Argentati M., Majorization for changes in angles between subspaces,
 {R}itz values, and graph {L}aplacian spectra, \href{https://doi.org/10.1137/060649070}{\textit{SIAM~J. Matrix Anal.
 Appl.}} \textbf{29} (2006), 15--32, \href{https://arxiv.org/abs/math.NA/0508591}{arXiv:math.NA/0508591}.

\bibitem{AJK2010}
Knyazev A., Jujunashvili A., Argentati M., Angles between infinite dimensional
 subspaces with applications to the {R}ayleigh--{R}itz and alternating
 projectors methods, \href{https://doi.org/10.1016/j.jfa.2010.05.018}{\textit{J.~Funct. Anal.}} \textbf{259} (2010), 1323--1345,
 \href{https://arxiv.org/abs/0705.1023}{arXiv:0705.1023}.

\bibitem{KT}
Kozlov V.V., Treshch\"ev D.V., Billiards. A genetic introduction to the
 dynamics of systems with impacts, \textit{Translations of Mathematical
 Monographs}, Vol.~89, \href{https://doi.org/10.1016/0021-8928(91)90103-2}{Amer. Math. Soc.}, Providence, RI, 1991.

\bibitem{M1}
Montgomery R., The three-body problem and the shape sphere, \href{https://doi.org/10.4169/amer.math.monthly.122.04.299}{\textit{Amer. Math. Monthly}} \textbf{122} (2015), 299--321, \href{https://arxiv.org/abs/1402.0841}{arXiv:1402.0841}.

\bibitem{MC1}
Murphy T.J., Cohen E.G.D., Maximum number of collisions among identical hard
 spheres, \href{https://doi.org/10.1007/BF01049961}{\textit{J.~Statist. Phys.}} \textbf{71} (1993), 1063--1080.

\bibitem{MC2}
Murphy T.J., Cohen E.G.D., On the sequences of collisions among hard spheres in
 infinite space, in Hard Ball Systems and the {L}orentz Gas,
 \href{https://doi.org/10.1007/978-3-662-04062-1_3}{\textit{Encyclopaedia Math. Sci.}}, Vol.~101, Springer, Berlin, 2000, 29--49.

\bibitem{Sev}
Sevryuk M.B., Estimate of the number of collisions of $n$ elastic particles on
 a~line, \href{https://doi.org/10.1007/BF01074110}{\textit{Theoret. and Math. Phys.}} \textbf{96} (1993), 818--826.

\bibitem{Sin}
Sinai Ya.G., Billiard trajectories in a polyhedral angle, \href{https://doi.org/10.1070/RM1978v033n01ABEH002254}{\textit{Russian Math.
 Surveys}} \textbf{33} (1978), 219--220.

\bibitem{Tab}
Tabachnikov S., Geometry and billiards, \textit{Student Mathematical Library},
 Vol.~30, \href{https://doi.org/10.1090/stml/030}{Amer. Math. Soc.}, Providence, RI, 2005.

\end{thebibliography}
\end{document}